# WAVELET THRESHOLDING FOR NONNECESSARILY GAUSSIAN NOISE: FUNCTIONALITY

BY R. AVERKAMP[1] AND C. HOUDRÉ[2]

*Freiburg University, and Université Paris XII and Georgia Institute of Technology*

For signals belonging to balls in smoothness classes and noise with enough moments, the asymptotic behavior of the minimax quadratic risk among soft-threshold estimates is investigated. In turn, these results, combined with a median filtering method, lead to asymptotics for denoising heavy tails via wavelet thresholding. Some further comparisons of wavelet thresholding and of kernel estimators are also briefly discussed.

**1. Introduction.** The model considered throughout these notes is the familiar one. The data takes the form

$$(1.1) \qquad X_i = f_i + \frac{e_i}{\sqrt{n}}, \qquad i = 1, \dots, n, n = 2^h, h \in \mathbb{N},$$

where $f = (f_i)$ is the signal to estimate and where the noise $e = (e_i)$ is such that the $e_i$ are zero mean i.i.d. random variables. One thinks of $f_i$ as $f_i = f_{n,i} = f(i/n)/\sqrt{n}$, so it is assumed that the data is sampled from a signal at the rate $1/n$ and then multiplied by $1/\sqrt{n}$. Applying a discrete wavelet transform (associated to an orthonormal wavelet basis, adapted to an interval and generated by a compactly supported wavelet) to the data leads to the noisy wavelet coefficients

$$(1.2) \qquad w_k = \theta_k + z_k; \qquad k = 1, \dots, 2^{j_0},$$

and

$$(1.3) \qquad w_{j,k} = \theta_{j,k} + z_{j,k}; \qquad j_0 \le j \le h-1, k = 1, \dots, 2^j,$$

Received September 2002; revised November 2004.

[1]Supported in part by NSF Grant DMS-96-32032.

[2]Supported in part by NSF Grant DMS-98-03239.

*AMS 2000 subject classifications.* Primary 62G07, 62C20; secondary 60G70, 41A25.

*Key words and phrases.* Wavelets, thresholding, minimax.

where to simplify notation the dependence on $n$ has been omitted (in particular, a factor $1/\sqrt{n}$ is omitted). Thresholding is then applied to the transformed data and the signal is recovered by applying an inverse transformation to the thresholded data [8]. In contrast to the ideal framework [5], in the functional framework the performance of estimators is no longer compared to a benchmark but instead the possible values of $(\theta_{\cdot,\cdot})$ is restricted to belonging to a ball in a smoothness class. To be more precise, it is assumed that

$$(1.4) \qquad \left(\sum_k |\theta_k|^p\right)^{1/p} + \left(\sum_{j \geq j_0}\left(2^{js}\left(\sum_k |\theta_{j,k}|^p\right)^{1/p}\right)^q\right)^{1/q} \leq A,$$

for some constant $A$, where $s := m + 1/2 - 1/p$ and $m > 1/p$. The condition $m > 1/p$ ensures that we deal with well-defined real-valued functions (and not generalized ones) in the Besov space $B_{p,q}^m$. [Recall that $m$ is the degree of smoothness of the function whose modulus of smoothness is locally quantified via (norms involving) the parameters $p$ and $q$.] Next, if the $(\theta_{\cdot,\cdot})$ are the wavelet coefficients of $f \in B_{p,q}^m$ and if the wavelet basis is sufficiently smooth, then $\|f\|_{B_{p,q}^m} \leq C_1 A$, where $C_1 = C_1(m, p, q)$ is a constant and $1 \leq p, q \leq +\infty$. Also, considering quasi-norms rather than norms, similar results hold in the cases $0 < p < 1$, or $0 < q < 1$ (we refer the reader to [3] and [15] for a much more extensive and precise list of references and further information on wavelets (and functions spaces)). We also note here that the Besov assumption can be replaced by a Triebel–Lizorkin one throughout much of the paper. Indeed, it is well known that the equivalence between the sequence space (quasi-)norm and the function space (quasi-)norm is what matters here. In view of this equivalence, we will slightly abuse notation and use $\| \cdot \|_{B_{p,q}^m}$ for the norm on the sequence space.

In this framework, Donoho, Hall, Johnstone, Kerkyacharian, Picard, Silverman and Yu compute minimax bounds of estimation [6, 7, 8, 9, 10, 11, 12], and show the corresponding optimality of wavelet thresholding. In particular, if the $e_i$ (hence, the $z_i$) are i.i.d. normal random variables, then the minimax rate in this model is $n^{-2m/(2m+1)}$, that is,

$$(1.5) \qquad \inf_{\hat{\theta}} \sup_{\theta \, : \, \|\theta\|_{B_{p,q}^m} \leq A} E\|\hat{\theta} - \theta\|_2^2 \sim C n^{-2m/(2m+1)},$$

where the infimum is taken over all estimators and where $C$ is a positive constant which depends on the variance of the noise, as well as on $m$, $p$, $q$ and $A$. [Throughout these notes, $\| \cdot \|_2$ is the Euclidean norm. From Parseval's identity and the equivalence between sequence and functional spaces, we thus see that (1.5) has an equivalent formulation at the function space level.] Moreover, estimators based on soft thresholding achieve this rate.

These early results were then extended to some classes of non-Gaussian noise by Neumann and Spokoiny [16] and Delyon and Juditsky [4]. It is shown in [16] that, for noise having finite moments of all orders (and $L^2$-differentiable density), soft thresholding achieves the same rate as soft thresholding for Gaussian noise. Furthermore, the actual performance, not just the rate, is the same. In [4], it is shown that, for more general distributions, soft thresholding can achieve the same rate as soft thresholding in the Gaussian case. In addition, under somewhat stronger conditions, the ratio of the minimax risk for Gaussian noise and other types of noise tends to one [16].

It is our purpose to further explore these topics here. Let us briefly discuss the contribution of the present paper. First, in Section 2 we show that if the noise only fulfills some moment conditions, soft thresholding actually achieves the same asymptotic performance as soft thresholding in the Gaussian case. In fact, it is shown that for soft thresholding the lim inf of the ratio of the minimax risk for Gaussian noise and this type of noise is larger than one. These results are then used, in Section 3, to tackle the estimation problem for noise with heavy tails. By first median filtering the data, the previous moment conditions become satisfied and then applying wavelet thresholding, it is still possible to have the same minimax rate as in the Gaussian case. To complete our study of wavelet thresholding methods, we return to the normal framework and present some concluding remarks comparing thresholding and kernel estimators with varying bandwidth.

**2. Moment conditions.** Our first statement is the core result of this section. To prove it, a fair amount of technical preparation is needed and the main part of the proof is postponed to the Appendix. However, we state and prove below some preparatory lemmas and indicate their use in the proof of the theorem.

In the sequel $\Phi$ denotes the standard normal distribution function, and $E_\Phi$ is expectation with respect to $\Phi$. Using the notation of [1], for any $\lambda > 0$, $T_\lambda^S$ denotes the soft thresholding operator given by $T_\lambda^S(x) = (|x| - \lambda)^+ \operatorname{sgn}(x), x \in \mathbb{R}$. Also, throughout the section the wavelet transform is as in [1], Section 4; in particular, the wavelet is assumed to be Hölder continuous of index $\beta > 0$.

THEOREM 2.1. *Let the model be given via* (1.1)–(1.4), *where* $p, q \geq 1$ *and* $m > 1/p$ *and where the* $e_i$ *have variance one. Let also the* $e_i$ *have finite moments of order* $L$, *where* $L$ *is such that*

$$(2.1) \qquad L > \frac{6}{2m/(2m+1)} \qquad if \ p \geq 2,$$

*and*

$$(2.2) \qquad L > \frac{6(m+1/2-1/p)(2m+1)}{(m+1/2-1/p)(2m+1)-m} \qquad if \ 1 \leq p \leq 2.$$

*Moreover, let the $e_i$ be symmetric. Then*

$$(2.3) \quad \liminf_{n \to \infty} \frac{\inf_{(\lambda) \in \mathbb{R}^n} \sup_{\theta : \|\theta\|_{B_{p,q}^m} \leq A} E_\Phi \sum_{j,k} (T_{\lambda_{j,k}}^S(w_{j,k}) - \theta_{j,k})^2}{\inf_{(\tilde{\lambda}) \in \mathbb{R}^n} \sup_{\theta : \|\theta\|_{B_{p,q}^m} \leq A} E \sum_{j,k} (T_{\tilde{\lambda}_{j,k}}^S(w_{j,k}) - \theta_{j,k})^2} \geq 1.$$

Above, the requirement of symmetry is imposed for technical reasons (we preserve the zero mean property of the wavelet transform of truncated noise). This requirement can be circumvented by more technical efforts in the proof. The i.i.d. assumption on the noise $e$ is not really needed either. Independence and $\sup_i E|e_i|^L < +\infty$, where $L$ satisfies either (2.1) or (2.2), will do, with also a variance level of $0 < \sigma^2 = \sup_i Ee_i^2 < +\infty$.

Let us illustrate the moment conditions to be satisfied: $p > 2 \Rightarrow L = 12$; $p = 2 \Rightarrow L > 12$; $p = 2, m \to +\infty \Rightarrow L > 6$; $p = 1 \Rightarrow L > 18$; $p = 1, m \to +\infty \Rightarrow L > 6$; $p = 3/2, m = 1 \Rightarrow L > 10$; $p = 3/2, m = 2 \Rightarrow L > 7.7$. Note that $L > 6$ is the least moment condition imposed above.

First, a well-known lemma whose proof is omitted.

LEMMA 2.2. *Let $s := m + 1/2 - 1/p$ and let*

$$\left( \sum_k |\theta_k|^p \right)^{1/p} + \left( \sum_{j \geq j_0} \left( 2^{js} \left( \sum_k |\theta_{j,k}|^p \right)^{1/p} \right)^q \right)^{1/q} \leq A,$$

*for some $A > 0$. Then for all $l \geq j_0$,*

$$\sum_{j \geq l} \|\theta_{j,\cdot}\|_2^2 \leq \begin{cases} A^2 (2^{-2m})^l / (1 - 2^{-2m}) = O(n^{-2\alpha m}), & \text{if } p \geq 2, \\ A^2 (2^{-2s})^l / (1 - 2^{-2s}) = O(n^{-2\alpha s}), & \text{if } 1 \leq p < 2, \end{cases}$$

*for any $\alpha$ such that $2^l \geq n^\alpha = 2^{\alpha h}$.*

As indicated in the Appendix, the previous lemma shows that, if we want to achieve the same minimax rate as in the Gaussian case, we need not worry about the (finer wavelet) coefficients in the levels $j \geq l = \alpha h$, as long as $\alpha > 1/(2m + 1)$, if $p \geq 2$, and $\alpha > m/((2m + 1)s)$ if $1 \leq p \leq 2$. Indeed, the square of the $\ell^2$-norm of these coefficients is of order $o(n^{-2m/(2m+1)})$. For $p \geq 2$, let $l$ be such that $2n^{1/(2m+1)} \geq 2^l > n^{1/(2m+1)}$. Then the simple estimator which discards the noisy coefficients of indices $l$ and above (keeping them otherwise) achieves the minimax rate since

$$\sum_{j \geq l, k} \theta_{j,k}^2 = O(n^{-2m/(2m+1)}) \quad \text{and} \quad \sum_{j < l, k} Ez_{j,k}^2 = O(n^{-2m/(2m+1)}).$$

Recall now a classical exponential inequality due to Kolmogorov (see [19], page 855).

LEMMA 2.3. *Let* $X_i$, $i = 1, \ldots, n$, *be zero mean, independent random variables. Let* $s_n^2 := \sum_{i=1}^n EX_i^2$, *let* $\sup_i \|X_i\|_\infty \leq K$, *and let* $S_n = \sum_{i=1}^n X_i$. *Then for all* $x > 0$,

$$P(S_n \geq s_n x) \leq \begin{cases} \exp\left(\dfrac{-x^2}{2}\left(1 - \dfrac{xK}{2s_n}\right)\right), & \text{if } x \leq s_n/K, \\ \exp\left(\dfrac{-xs_n}{4K}\right), & \text{if } x \geq s_n/K. \end{cases}$$

The next lemma is a simple application of the previous one. It is used in the proof of Theorem 2.1 to upper estimate $E(T^S_{\lambda_{j,k}}(z_{j,k} + \theta_{j,k}) - \theta_{j,k})^2 \mathbf{1}_{\{|z_{j,k}| > b_{j,k}\}}$, for appropriately chosen $\lambda_{j,k}$ and $b_{j,k}$.

LEMMA 2.4. *Let* $(X_{i,n})_{i,n\in\mathbb{N}}$ *be zero mean random variables such that, for each fixed* $n$, *the* $X_{i,n}$ *are independent. Let* $\sum_i EX_{i,n}^2 = 1$ *and let* $\sup_i \|X_{i,n}\|_\infty \leq K_n$, *where* $\lim_{n\to\infty} K_n = 0$. *Let* $F_n$ *be the distribution function of* $\sum_i X_{i,n}$, *and let* $(a_n)$ *be a sequence of positive reals with* $a_n = o(1/K_n)$ *and such that, for all* $n \in \mathbb{N}$, $k_n := (1 - a_n K_n/2) > 0$. *Then, for any* $a$ *with* $0 < a < a_n$,

$$\int_a^\infty x^2 F_n(dx) \leq \frac{a^2 + 2}{k_n} \exp(-k_n a^2/2) + o(\exp(-1/K_n)).$$

PROOF. Using Lemma 2.3, we have

$$\int_a^\infty x^2 \, dF_n(x) = a^2(1 - F_n(a)) + 2\int_a^\infty x(1 - F_n(x)) \, dx$$

$$\leq a^2 \exp(-k_n a^2/2) + 2\int_a^\infty x \exp(-k_n x^2/2) \, dx$$

$$+ 2\int_{1/K_n}^\infty x \exp(-x/(4K_n)) \, dx$$

$$= a^2 \exp(-k_n a^2/2) + 2/k_n \exp(-k_n a^2/2)$$

$$- [8K_n x \exp(-x/(4K_n))]_{1/K_n}^\infty + 8K_n \int_{1/K_n}^\infty \exp(-x/(4K_n)) \, dx$$

$$= a^2 \exp(-k_n a^2/2) + 2/k_n \exp(-k_n a^2/2)$$

$$+ 8\exp(-1/(4K_n^2)) + 32K_n^2 \exp(-1/(4K_n^2))$$

$$\leq (a^2 + 2)/k_n \exp(-k_n a^2/2) + o(\exp(-1/K_n)). \qquad \square$$

We further need the following large deviation result, which is a simple extension of Lemma 5.8 in [17]; the difference with this lemma is that the requirement of identical distributions is dropped. The proof with the help of Esseen's inequality ([17], Theorem 5.4) is essentially the same as for Lemma 5.8 in [17] (our $C$ below is $A$ in [17]).

This lemma is used to show that, for a large class of noise, and midsize thresholds, the soft thresholding risk converges to the Gaussian risk.

LEMMA 2.5. *Let* $(X_{i,n})_{i,n \in \mathbb{N}}$ *be zero mean random variables such that, for each fixed* $n$, *the* $X_{i,n}$ *are independent. Let* $\sum_i EX_{i,n}^2 = 1$ *and let* $M_n := \sum_i E|X_{i,n}|^3 < +\infty$. *Then for all* $0 < \varepsilon < 1$ *there exist* $\beta_n$ *with* $\beta_n \to 1$ *such that, for all* $x$ *with* $|x| \le (1 - \varepsilon)\sqrt{2 \log(1/(CM_n))}$,

$$(2.4)\quad \beta_n \le \frac{P(\sum_i X_{i,n} \le x)}{\Phi((-\infty, x])} \le 1/\beta_n \quad \text{and} \quad \beta_n \le \frac{P(\sum_i X_{i,n} > x)}{\Phi((x, +\infty))} \le 1/\beta_n,$$

*where* $C$ *is an absolute constant.*

REMARK 2.6. From the end of the proof of Theorem 2.1, we infer that the thresholds $\tilde{\lambda}_{j,k}$ can asymptotically be chosen as in the Gaussian case ($\tilde{\lambda} = \lambda$ for $\lambda < \alpha_n/2$, while $\lambda/\tilde{\lambda} \to 1$ for $\lambda \ge \alpha_n/2$). However, in this proof the thresholds $\tilde{\lambda}$ are larger than the Gaussian ones (with the same variance). In the Gaussian functional approach, the optimal minimax thresholds are of order $C\sigma\sqrt{(j - j_0)_+}/\sqrt{n}$, where $C$ and $j_0$ depend on $m, p$ and $q$ and where $\sigma^2$ is the variance of the noise (e.g., see [4, 16]). For the ideal estimator approach, the optimal minimax rate is achieved with thresholds of uniform size $\sim \sigma\sqrt{2 \log n}/\sqrt{n}$, and we also know (see Theorem 6.1 in [1]) that thresholds can be chosen levelwise to still produce a minimax method. There, for the level $j$ the thresholds were chosen to be of size $\sim \sigma\sqrt{2j \log 2}/\sqrt{n}$. Now, using thresholds of size $C\sigma\sqrt{j}/\sqrt{n}$ for the level $j$ in the function space approach almost achieves the ideal minimax rate; it is only worse by a factor $O(\log n)$. This discrepancy cannot be avoided in general and, at least for $p \ge 2$, no set of thresholds will achieve the optimal minimax rate in both contexts. Indeed, let $X_i = f_i + e_i$, $i = 1, \ldots, n$, where the $e_i$ are i.i.d. normal random variables with mean zero and variance $1/n$, and let $p_\Phi(\cdot, \cdot)$ be defined as in the Appendix (or as Theorem 2.1 in [1]). If $\lambda \ge \theta \ge 0$, then

$$(2.5)\quad p_\Phi(\lambda, \theta) \ge \theta^2 \Phi((-\lambda - \theta, \lambda - \theta)) + \int_{-\infty}^{-\lambda - \theta}(x + \lambda)^2 \Phi(dx) \ge \frac{\theta^2}{2}.$$

Let $\lambda_{n,j}$ ($n$ is for the number of coefficients, while $j$ is a particular level) be a set of thresholds which achieve the optimal minimax rate in the ideal estimator context. For a fixed $\alpha \in (0, 1)$, the optimal thresholds for the level $j = \alpha \log_2 n$ have to be at least of size $\sim \sqrt{2j \log 2}/\sqrt{n} = C\sqrt{j}/\sqrt{n}$, where $C$ is a constant. The reason is that $2^j p_\Phi(\lambda_j, 0) = O(\log n/n)$ is needed to achieve the minimax rate for the ideal estimator approach. Let now

$$j_0 := \min_j \left\{ j : \frac{C\sqrt{j}}{\sqrt{n}} \ge 2A\sqrt{2^{-j(2m+1)}} \right\}.$$

Simple computations yield that $j_0 \sim (\log_2 n)/(2m+1)$. If $\theta_{j_0,k} = A\sqrt{2^{-j_0(2m+1)}}$, $k = 0, \ldots, 2^{j_0} - 1$, and $\theta_{j,k} = 0$ elsewhere, then clearly $\|\theta\|_{B_{p,q}^m} \leq A$. If $n$ tends to infinity, then for $n$ larger than a certain bound, $\lambda_{n,j_0} > A\sqrt{2^{-j_0(2m+1)}}$. Now, it follows from (2.5) that the risk for thresholding the signal $(\theta.)$ at level $j_0$ with thresholds $\lambda_{n,j_0}$ is at least as large as

$$2^{j_0} A^2 2^{-j_0(2m+1)}/2 = A^2 2^{-j_0 2m}/2. \tag{2.6}$$

Using the definition of $j_0$, we obtain

$$2^{2m+3} A^2 2^{-j_0(2m+1)} \geq \frac{C^2(j_0 - 1)}{n}. \tag{2.7}$$

Combining the relations (2.6) and (2.7) shows that the risk for estimating $(\theta.)$ is as large as

$$A^2 2^{-j_0 2m - 1}$$
$$\geq A^2 \left(\frac{C}{A}\right)^{4m/(2m+1)} n^{-2m/(2m+1)} (j_0 - 1)^{2m/(2m+1)} 2^{-2m(2m+3)/(2m+1) - 1}.$$

Since $j_0 \sim \log_2 n/(2m+1)$, this is worse than the minimax rate for $B_{p,q}^m$.

## 3. Heavy tails and median filtering.

To date, asymptotics for wavelet thresholding seems to have been restricted to noise with higher moments. Next, we want to try to apply wavelet thresholding to noise with heavy tails and study the corresponding quadratic risks. By first applying a median filter to the data, the absence of finite moments will be overcome. The downside of this approach, however, is that it introduces an additional bias. Nevertheless, under these conditions wavelet thresholding applied to the filtered data achieves at least the same minimax rate as in the normal case, but the constants are larger. Various types of nonlinear smoothers involving medians have proved useful in time series analysis (e.g., see [14, 18, 20]). Another, wavelet inspired, approach to denoising heavy tails based on a different preprocessing method is also developed in [10].

Below, given $a_1, \ldots, a_{2k+1}$, let $\text{med}(a_1, \ldots, a_{2k+1})$ be the real $x$ such that $\#\{i : a_i \geq x\} = k+1$ and $\#\{i : a_i \leq x\} = k+1$, with $\#$ denoting cardinality. To simplify notation, we use the abbreviation $\text{med}(a_i, 2k+1)$ for $\text{med}(a_{i-k}, \ldots, a_{i+k})$. If $i$ is smaller than $k$, then $\text{med}(a_i, 2k+1) := \text{med}(a_1, \ldots, a_{2k+1})$; a similar boundary correction is performed for the largest indices.

Our first lemma makes the advantage of the median filter clear as far as the existence of moments is concerned. It shows, for example, that the median of thirteen independent Cauchy random variables has moments of order $7 - \varepsilon$, $\varepsilon > 0$.

LEMMA 3.1. *Let $X_1, \ldots, X_{2k-1}$, $k \geq 1$, be independent random variables. For any $x > 0$,*

$$P(\mathrm{med}(X_1, \ldots, X_{2k-1}) \geq x) \leq \binom{2k-1}{k} \max_{i=1,\ldots,2k-1} (P(X_i \geq x))^k.$$

*In particular, if there exist constants $C > 0, \gamma > 0$ such that, for $x$ large enough, $\max_{i=1,\ldots,2k-1} P(|X_i| \geq x) \leq \frac{C}{x^\gamma}$, then $\mathrm{med}(X_1, \ldots, X_{2k-1})$ has moments of order $r < k\gamma$.*

PROOF.

$$\{\mathrm{med}(X_1, \ldots, X_{2k-1}) \geq x\} = \bigcup_{\substack{M \subset \{1,\ldots,2k-1\} \\ \#M=k}} \{X_i \geq x : i \in M\}.$$

Hence,

$$P(\mathrm{med}(X_1, \ldots, X_{2k-1}) \geq x) \leq \binom{2k-1}{k} \max_{i=1,\ldots,2k-1} (P(X_i \geq x))^k. \quad \square$$

Let us now give the main result of this section. As before, the data is given via (1.1), while $W_n$ is a discrete wavelet transform as in the previous section (in particular, it is generated by a compactly supported wavelet $\psi$ which is Hölder continuous of index $\beta > 0$, chosen later). Again, let $\theta = W_n(f)$ and let also $\mathrm{med}(X, 2l+1) := (\mathrm{med}(X_i, 2l+1))_{1 \leq i \leq n}$, where (using the notation set above) $\mathrm{med}(X_i, 2l+1) := \mathrm{med}(X_{i-l}, \ldots, X_{i+l})$, for $i - l \geq 1$ and $\mathrm{med}(X_i, 2l+1) := \mathrm{med}(X_1, \ldots, X_{2l+1})$ otherwise.

THEOREM 3.2. *Let the $e_i$ be symmetric with $E|e_1|^\gamma < +\infty$, for some $\gamma > 0$. Let $A, B > 0$. Then there exist an $l = l(\gamma)$ and thresholds $\lambda_{j,k}$ such that*

$$\sup_{\substack{\|\theta\|_{B_{p,q}^m} \leq A \\ \sum_i |f_i - f_{i-1}|^2 \leq B/n}} E \sum_{j,k} |T_{\lambda_{j,k}}^S (W_n(\mathrm{med}(X, 2l+1))_{j,k}) - \theta_{j,k}|^2 = O(n^{-2m/(2m+1)}).$$

We impose the condition $\sum_i |f_i - f_{i-1}|^2 \leq B/n$ to have control over the $\ell_2$-norm of the bias, that is, on

$$(3.1) \qquad \sum_{i=1}^n (\mathrm{med}(X_i, 2l+1) - \mathrm{med}(e_i, 2l+1) - f_i)^2,$$

which we introduce by median filtering the data. This condition is not that strong and in most cases follows from the Besov norm condition. We will take another look at this after the proof of the theorem.

Theorem [3.2] is more than just an existence result. Indeed, from its proof we infer that the above thresholds $\lambda_{j,k}$ can asymptotically be chosen as in the Gaussian case, but with a new variance which is now at most $2D^2\sigma_{\max}^2$, with $\sigma_{\max}^2$ given below and with $D = 2l + 1$ (see also Remark [2.6]).

PROOF OF THEOREM [3.2]. Since the $e_i$ are symmetric, $E\operatorname{med}(e_i, 2l+1) = 0$ (again $l$ is chosen later). Let $y_{j,k}$ be the coefficient of index $j, k$ of $W_n(\operatorname{med}(e, 2l+1))$, and let

$$(b_{j,k}) := W_n(\operatorname{med}(X, 2l+1)) - W_n(f) - (y_{j,k}).$$

First we prove that the influence of the random variables $b_{j,k}$ (the bias) is not too large in our estimation problem:

$$E(T^S_{\lambda_{j,k}}(W_n(\operatorname{med}(X, 2l+1))_{j,k}) - \theta_{j,k})^2$$

$$= E(T^S_{\lambda_{j,k}}(\theta_{j,k} + b_{j,k} + y_{j,k}) - \theta_{j,k})^2$$

$$= E(T^S_{\lambda_{j,k}}(\theta_{j,k} + b_{j,k} + y_{j,k}) - T^S_{\lambda_{j,k}}(\theta_{j,k} + y_{j,k}) + T^S_{\lambda_{j,k}}(\theta_{j,k} + y_{j,k}) - \theta_{j,k})^2$$

$$\leq 2Eb_{j,k}^2 + 2E(T^S_{\lambda_{j,k}}(\theta_{j,k} + y_{j,k}) - \theta_{j,k})^2,$$

since $|T^S_\lambda(x_1) - T^S_\lambda(x_2)| \leq |x_1 - x_2|$. Thus,

$$\sum_{j,k} E(T^S_{\lambda_{j,k}}(\theta_{j,k} + b_{j,k} + y_{j,k}) - \theta_{j,k})^2$$

$$\leq 2\sum_{j,k} Eb_{j,k}^2 + 2\sum_{j,k} E(T^S_{\lambda_{j,k}}(\theta_{j,k} + y_{j,k}) - \theta_{j,k})^2.$$

Note that

$$\sum_{j,k} b_{j,k}^2 = \sum_{j,k} (W_n(\operatorname{med}(X, 2l+1) - \operatorname{med}(e, 2l+1) - f))_{j,k}^2$$

$$= \sum_{i=1}^n (\operatorname{med}(X_i, 2l+1) - \operatorname{med}(e_i, 2l+1) - f_i)^2.$$

But for $l < i \leq n - l$,

$$|\operatorname{med}(X_i, 2l+1) - \operatorname{med}(e_i, 2l+1) - f_i|$$

$$\leq |\operatorname{med}(e_{i-l} + f_{i-l} - f_i, \ldots, e_{i+l} + f_{i+l} - f_i) - \operatorname{med}(e_i, 2l+1)|$$

$$\leq \max_{j=-l,\ldots,l} |f_{i+j} - f_i|$$

$$\leq \sum_{j=-l+1}^l |f_{i+j} - f_{i+j-1}|.$$

If $i \leq l$ or $i > n - l$, then, similarly, $|\operatorname{med}(X_i, 2l+1) - \operatorname{med}(e_i, 2l+1) - f_i| \leq \sum_{j=2}^{2l+1} |f_j - f_{j-1}|$, respectively, $\leq \sum_{j=n-2l+1}^{n} |f_j - f_{j-1}|$. Hence,

$$
\begin{aligned}
\sum_{j,k} b_{j,k}^2 &\leq \sum_{i=l+1}^{n-l} 2l \sum_{j=-l+1}^{l} |f_{i+j} - f_{i+j-1}|^2 \\
&\quad + 2l^2 \sum_{j=2}^{2l+1} |f_j - f_{j-1}|^2 + 2l^2 \sum_{j=n-2l+1}^{n} |f_j - f_{j-1}|^2 \\
&\leq 8l^2 \sum_{i=2}^{n} |f_i - f_{i-1}|^2 = O(1/n).
\end{aligned}
$$

This implies that, if we choose a fixed median filter, then $\sum_{j,k} E b_{j,k}^2$ is negligible compared to $O(n^{-2m/(2m+1)})$. Thus, to finish the proof, it suffices to show that

$$
\sup_{\|\theta\|_{B_{p,q}^m} \leq A} E \sum_{j,k} (T_{\lambda_{j,k}}^S (W_n(f + \tilde{e})_{j,k}) - \theta_{j,k})^2 = O(n^{-2m/(2m+1)}),
$$

where $\tilde{e}_i := \operatorname{med}(e_i, D)$ and $D = 2l + 1$ is chosen such that $E|\tilde{e}_1|^L < \infty$ and $L$ satisfies the moment conditions (which depend on $\gamma$) of Theorem 2.1. If the $\tilde{e}_i$ were independent, which they are not, we could apply Theorem 2.1 to conclude. The next two lemmas deal with this new situation (the $D$-dependent case) and, respectively, correspond to Lemma 2.3 and to Lemma 2.5 in the independent case.

LEMMA 3.3. *Let $X_1, \ldots, X_n$ be zero mean bounded random variables, with $\sup_i \|X_i\|_\infty \leq K$, and also $D$-dependent, that is, such that $X_{i_1}, \ldots, X_{i_k}$ are independent if $\min_{1 \leq j \neq r \leq k} |i_j - i_r| \geq D$. Let $S_j = \sum_{i=0}^{[(n-1)/D]} X_{iD+j}$, $j = 1, \ldots, D$, $\sigma_j^2 = E S_j^2$, and $\sigma_{\max} = \max_{j=1,\ldots,D} \sigma_j$. Then*

$$
(3.2) \quad P\left(\sum_{i=1}^{n} X_i \geq x\right) \leq D
\begin{cases}
\exp\left(\dfrac{-x^2}{4D^2 \sigma_{\max}^2}\right), & \text{if } x \leq \dfrac{\sigma_{\max}^2 D}{K}, \\[2ex]
\exp\left(\dfrac{-x}{4KD}\right), & \text{if } x \geq \dfrac{\sigma_{\max}^2 D}{K} .
\end{cases}
$$

PROOF. Note that the $S_j$ are sums of independent random variables:

$$
\begin{aligned}
P\left(\sum_{i=1}^{n} X_i \geq x\right) &\leq \sum_{i=1}^{D} P(S_i \geq x/D) \\
&= \sum_{i=1}^{D} P(S_i/\sigma_i \geq x/(\sigma_i D))
\end{aligned}
$$

$$\leq \sum_{i=1}^{D} \begin{cases} \exp\left(\dfrac{-x^2}{2D^2\sigma_i^2}\left(1 - \dfrac{xK}{2D\sigma_i^2}\right)\right), & \text{if } x \leq \dfrac{\sigma_i^2 D}{K}, \\[3mm] \exp\left(\dfrac{-x}{4KD}\right), & \text{if } x \geq \dfrac{\sigma_i^2 D}{K}, \end{cases}$$

$$\leq \sum_{i=1}^{D} \begin{cases} \exp\left(\dfrac{-x^2}{2D^2\sigma_i^2}\dfrac{1}{2}\right), & \text{if } x \leq \dfrac{\sigma_i^2 D}{K}, \\[3mm] \exp\left(\dfrac{-x}{4KD}\right), & \text{if } x \geq \dfrac{\sigma_i^2 D}{K}, \end{cases}$$

$$\leq D \begin{cases} \exp\left(\dfrac{-x^2}{4D^2\sigma_{\max}^2}\right), & \text{if } x \leq \dfrac{\sigma_{\max}^2 D}{K}, \\[3mm] \exp\left(\dfrac{-x}{4KD}\right), & \text{if } x \geq \dfrac{\sigma_{\max}^2 D}{K}, \end{cases}$$

where the last inequality holds since for $x \leq \frac{\sigma_{\max}^2 D}{K}$, $\frac{x^2}{4D^2\sigma_{\max}^2} \leq \frac{x}{4KD}$. □

With the help of Lemma 3.3, it is also possible to prove a $D$-dependent version of Lemma 2.4.

Lemma 3.4. *Let $(X_{i,n})_{i,n\in\mathbb{N}}$ be zero mean random variables such that, for each fixed $n$, the $X_{i,n}$ are $D$-dependent and such that $M_n := \sum_i E|X_{i,n}|^3 < +\infty$. Let $S_{j,n} = \sum_i X_{iD+j,n}$, $j = 1,\ldots,D$, let $\sigma_{j,n}^2 = ES_{j,n}^2$ and let $\sigma_{\max,n} = \max_{j=1,\ldots,D} \sigma_{j,n}$. Then for all $0 < \varepsilon < 1$, there exist $\beta_n = \beta_n(\varepsilon)$ with $\limsup_{n\to+\infty} \beta_n = D$ such that*

$$\sup_{0\leq x\leq \varepsilon D\sigma_{\max,n}\sqrt{2\log(1/M_n)}} \frac{P(\sum_i X_{i,n} \leq -x)}{\Phi((-\infty, -x/(\sigma_{\max,n}D)))} \leq \beta_n$$

*and*

$$\sup_{0\leq x\leq \varepsilon D\sigma_{\max,n}\sqrt{2\log(1/M_n)}} \frac{P(\sum_i X_{i,n} \geq x)}{\Phi((x/(\sigma_{\max,n}D), +\infty))} \leq \beta_n.$$

Proof. If $x \geq 0$, then

$$P\left(\sum_i X_{i,n} \geq x\right) \leq \sum_{j=1}^{D} P(S_{j,n} > x/D)$$

$$\leq \sum_{j=1}^{D} P(S_{j,n}/\sigma_{j,n} > x/(\sigma_{\max,n}D)).$$

A similar inequality holds for $x \leq 0$. Since the $S_{j,n}$ are sums of independent random variables, the assertion follows from Lemma 2.5. □

Note that it is, moreover, trivial that, for $x \geq 0$,

$$\frac{P(\sum_i X_{i,n} \leq x)}{\Phi((-\infty, x/(\sigma_{\max,n}D)))} \leq 2 \quad \text{and} \quad \frac{P(\sum_i X_{i,n} \geq -x)}{\Phi((-x/(\sigma_{\max,n}D), +\infty))} \leq 2;$$

and this gives a version of the other half of Lemma 2.5 in the $D$-dependent case.

The rest of the proof of Theorem 3.2 is then quite similar to the proof of Theorem 2.1. Let us return to it. Again, let $\tilde{e}_i := \operatorname{med}(e_i, D)$ be as defined above. First, as in the proof of Theorem 2.1, we can assume that the $\tilde{e}_i$ are bounded by $n^\delta$ for some $\delta > 0$, since the upper estimate in Lemma A.1 holds for $D$-dependent random variables with a constant depending now also on $D$.

Of importance in the proof of Theorem 2.1 was the distribution of the noise in the wavelet coefficients. We denote the coefficients of the wavelet transform by $(c_{j,k,i})$, that is, $\theta_{j,k} = \sum_i c_{j,k,i} f_i$. At the boundary we have the problem that $\tilde{e}_1 = \cdots = \tilde{e}_{(D+1)/2}$ and $\tilde{e}_{n-(D-1)/2} = \cdots = \tilde{e}_n$. But

$$y_{j,k} = \left(\tilde{e}_1 \sum_{i=1}^{(D+1)/2} c_{j,k,i}\right) + \sum_{i=(D+1)/2+1}^{n-(D-1)/2-1} c_{j,k,i}\tilde{e}_i + \left(\tilde{e}_n \sum_{i=n-(D-1)/2}^{n} c_{j,k,i}\right),$$

and this last expression, which is a sum of $n - D + 1$ random variables which are $D$-dependent, thus satisfies (after reordering) the conditions of Lemmas 3.3 and 3.4. Anyway, only about $O(\log n)$ wavelet coefficients are affected by this problem. If we do not threshold these coefficients, the risk would increase at most by $O((\log n)/n)$ and this is negligible compared to the minimax risk. Let us using $y_{j,k,r} = \sum_i c_{j,k,iD+r}\tilde{e}_{iD+r}$. Then using the $D$-dependence condition, we see that $y_{j,k,r}$ is a sum of independent random variables with, moreover, $y_{j,k} = \sum_{r=1}^{D} y_{j,k,r}$. Let $\sigma_{j,k,r}^2 = E y_{j,k,r}^2$ and let $\sigma_{j,k,\max}^2 = \max_r \sigma_{j,k,r}^2$. Then, clearly, $\sigma_{j,k,\max}^2 \leq E\tilde{e}_1^2$ and a version of Lemma A.1 holds for the $y_{j,k}$, the upper constants depending now also on $D$. Given this, as well as Lemmas 3.3 and 3.4, we can now proceed as in the proof of Theorem 2.1. Hence, with the right thresholds, we can achieve $D$ times the performance of the Gaussian risk with variance $2\sigma_{\max}^2 D^2$, where $\sigma_{\max}^2 = \max_{j \leq t,k} \sigma_{j,k,\max}^2$, $t$ being the finest level where the wavelet coefficient is not discarded.

Also of importance in Lemmas 3.3 and 3.4 is the term $\sigma_{\max}^2$; we show next that in general,

$$(3.3) \qquad\qquad \sigma_{j,k,\max}^2 \approx E\tilde{e}_1^2/D.$$

Let $h = \log_2 n$. Since the wavelet $\psi$ is compactly supported and Hölder continuous of index $\beta$, we know that (see, e.g., the proof of Theorem 4.1 in [1])

$$|2^{(h-j)/2}c_{j,k,i} - \psi(2^{j-h}i - k)| \leq C_1 2^{(j-h)\beta},$$

with also

$$|\psi(2^{j-h}i - k) - \psi(2^{j-h}(i-1) - k)| \leq C_2 2^{(j-h)\beta},$$

for some constants $C_1, C_2$. Thus,

$$|c_{j,k,i} - c_{j,k,i+1}| \leq (2C_1 + C_2)2^{(j-h)(\beta+1/2)}$$

and

$$|c_{j,k,i}^2 - c_{j,k,i+1}^2| = |c_{j,k,i} - c_{j,k,i+1}||c_{j,k,i} + c_{j,k,i+1}|$$
$$\leq C_3 2^{(j-h)/2} 2^{(j-h)(\beta+1/2)},$$

since the $|c_{j,k,i}|$ are of order $O(2^{(j-h)/2})$, $C_3$ being a constant. But $\sigma_{j,k,r}^2 = E\tilde{e}_1^2 \sum_i c_{j,k,iD+r}^2$ and $\#\{i : c_{j,k,i} \neq 0\} = O(2^{h-j})$ (again, $\#$ denotes cardinality) since the wavelet is compactly supported. Indeed, recall that (see, e.g., [3])

$$(3.4) \qquad \phi_{j,k} = \sum_{i=0}^{2^{j_0-j}(N-1)} u_{j_0-j, i+2^{j_0-j}k}\phi_{j_0,i}$$

and

$$(3.5) \qquad \psi_{j,k} = \sum_{i=0}^{2^{j_0-j}(N-1)} v_{j_0-j, i+2^{j_0-j}k}\phi_{j_0,i},$$

where $u_{\cdot,\cdot}$ and $v_{\cdot,\cdot}$ depend only on the scaling identities (whose size we set equal to $N$). This claim about the length of the filters $(u_{j_0-j})$ and $(u_{j_0-j})$ can be proved via a simple induction argument. Actually (see [1]), $\max_i |u_{j_0-j,i}| = O(2^{(j_0-j)/2})$ and $\max_i |v_{j_0-j,i}| = O(2^{(j_0-j)/2})$. Thus,

$$|\sigma_{j,k,r}^2 - \sigma_{j,k,r+1}^2| = \left| E\tilde{e}_1^2 \sum_i (c_{j,k,iD+r}^2 - c_{j,k,iD+r+1}^2) \right|$$
$$= O(2^{(j-h)\beta}|E\tilde{e}_1^2|).$$

Since $\sum_r \sigma_{j,k,r}^2 = E\tilde{e}_1^2 \sum_i c_{j,k,i}^2 = E\tilde{e}_1^2$, all the $\sigma_{j,k,r}^2$ have about the same size and, thus, $\sigma_{j,k,\max}^2 \approx E\tilde{e}_1^2/D$. This completes the proof of Theorem 3.2. $\square$

We now turn to the problem of finding out when the condition

$$(3.6) \qquad \sum_{i=1}^{n-1} |f_i - f_{i+1}|^2 = O(1/n)$$

follows from $\|\theta\|_{B_{p,q}^m} \leq A$. If $m \leq 1$, then assume $\beta \geq m$, where $\beta \leq 1$ is the Hölder continuity exponent of the wavelet (otherwise the characterization of smoothness via wavelets does not make sense). Again, $h = \log_2 n$. Since for a constant $C_1 > 0$,

$$|c_{j,k,i} - c_{j,k,i+1}| \leq C_1 2^{(j-h)(1/2+\beta)},$$

and $\#\{i : c_{j,k,i} \neq 0\} = O(2^{h-j})$, it follows that

$$\sum_{i=1}^{n-1} |c_{j,k,i} - c_{j,k,i+1}|^2 \leq C_2 2^{2\beta(j-h)},$$

where $C_2$ is another constant. Note that since the wavelet transform is an orthonormal transformation, $f_i = \sum_{j,k} a_{j,k} c_{j,k,i}$. Thus,

$$\begin{aligned}
\sum_{i=1}^{n-1} |f_i - f_{i+1}|^2 &= \sum_{i=1}^{n-1} \left( \sum_{j,k} a_{j,k}(c_{j,k,i} - c_{j,k,i+1}) \right)^2 \\
&\leq \sum_{i=1}^{n-1} h \sum_{j=0}^{h-1} \left( \sum_k a_{j,k}(c_{j,k,i} - c_{j,k,i+1}) \right)^2 \\
&\leq \sum_{i=1}^{n-1} h \sum_{j=0}^{h-1} C_3 \sum_k (a_{j,k}(c_{j,k,i} - c_{j,k,i+1}))^2,
\end{aligned}$$

since $\#\{k : c_{j,k,i} \neq 0 \text{ or } c_{j,k,i+1} \neq 0\} = O(1)$, see (3.4) and (3.5), and with $C_3$ a constant

$$\begin{aligned}
&= hC_3 \sum_{j=0}^{h-1} \sum_k a_{j,k}^2 \sum_i (c_{j,k,i} - c_{j,k,i+1})^2 \\
&\leq hC_3 C_2 \sum_{j=0}^{h-1} 2^{2\beta(j-h)} \sum_k a_{j,k}^2 \\
&\leq hC_3 C_2 A^2 \sum_{j=0}^{h-1} 2^{2\beta(j-h)} \begin{cases} 2^{-2jm}, & \text{if } p \geq 2, \\ 2^{-2js}, & \text{if } p \leq 2, \end{cases}
\end{aligned}$$

where the last inequality is proved by using arguments as in the proof of Lemma 2.2. Thus, if $\beta = 1$ and $m \geq 1$, respectively, $s \geq 1$, then the last term is equal to $O(hn^{-2})$. If $m \leq \beta$, respectively, $s \leq \beta$, then the last term is equal to $O(hn^{-2m})$, respectively, $O(hn^{-2s})$. Hence, for $p \geq 2$ we obtain

$$\sum_{i=1}^{n-1} |f_i - f_{i+1}|^2 = O(\log n / n^{(2m \wedge 2)}), \tag{3.7}$$

and for $p \leq 2$ we obtain

$$\sum_{i=1}^{n-1} |f_i - f_{i+1}|^2 = O(\log n / n^{(2s \wedge 2)}). \tag{3.8}$$

Thus, for $p \geq 2$ the condition $m > 1/2$ will ensure that (3.6) holds. For $p \leq 2$ the additional condition $m > 1/p$ ensures that $2s > 1$ and, thus, (3.6) is always satisfied.

REMARK 3.5.   Above, and also in view of the proof of Theorem 2.1, the i.i.d. assumption on $e$ can be weakened and replaced by independence with $\sup_i E|e_i|^\gamma < +\infty$, for some $\gamma > 0$. The previous proofs also show how to deal with noise (with or without higher moments) that is not independent, but $D$-dependent, where $D$ is a fixed constant. Indeed, Lemmas 3.3 and 3.4 are applicable, and then it is easy to mimic the proof of Theorem 3.2 and the minimax rate for this problem is again as in the Gaussian case. To obtain such a result, the noisy wavelet coefficients need not converge in distribution to a normal random variable. Only the bounds of Lemma 3.4 and of Lemma 3.3 are needed. This approach via large deviation results is also possible for other kinds of correlated noise. Under appropriate weak dependence conditions, the law of the noisy wavelet coefficients is asymptotically normal with a variance possibly bigger than the variance of the original noise. Wavelet thresholding has also been investigated for stationary Gaussian noise; for example, see [12, 21]. Let us finally mention that it would be interesting to transfer the "ideal framework with quadratic risk" to heavy tail noise via median filtering.

REMARK 3.6.   The upper bounds obtained in Theorems 2.1 and 3.2 can often be complemented with lower bounds of the same order for various types of noises. In turn, these bounds often represent the order of the minimax rate among all estimators (see the various references cited in the introductory section). However, different nonlinear estimators can outperform wavelet thresholding for still other types of noise. Let us briefly present such an estimator. The model is the usual one, $X_i = f_i + e_i/\sqrt{n}$, $i = 1, \ldots, n = 2^h$, where the $e_i$ are zero mean i.i.d. random variables with finite second moment. Our estimator of $f_i$ based on the $X_i$ is

$$(3.9) \quad \hat{f}_i := \max_{j=0,\ldots,M-1} X_{i+j} - \frac{E \max_{i=1,\ldots,M} e_i}{\sqrt{n}}, \qquad i = 1, \ldots, n - M + 1,$$

and for $i > n - M + 1$,

$$(3.10) \qquad\qquad\qquad \hat{f}_i = \hat{f}_{n-M+1},$$

where $M := M(n)$ will be chosen later. Let $c_M := E \max_{i=1,\ldots,M} e_i$. Thus, for $i \leq n - M + 1$,

$$\hat{f}_i - f_i = \max_{j=0,\ldots,M-1} \left( f_{i+j} - f_i + \frac{1}{\sqrt{n}}(e_{i+j} - c_M) \right).$$

Hence,

$$|\hat{f}_i - f_i| \le \max_{j=0,\dots,M-1} |f_{i+j} - f_i| + \frac{1}{\sqrt{n}} \Big| \max_{j=0,\dots,M-1} e_{i+j} - c_M \Big|$$

and

$$E|\hat{f}_i - f_i|^2 \le 2\left(\sum_{j=1}^{M-1} |f_{i+j} - f_{i+j-1}|\right)^2 + \frac{2}{n} E\left(\max_{j=0,\dots,M-1} e_{i+j} - c_M\right)^2$$

$$\le 2M \sum_{j=1}^{M-1} |f_{i+j} - f_{i+j-1}|^2 + \frac{2}{n} E\left(\max_{j=0,\dots,M-1} e_{i+j} - c_M\right)^2.$$

A similar computation for $i > n - M + 1$ gives

$$E|\hat{f}_i - f_i|^2 \le 2M \sum_{j=m-M+1}^{n} |f_j - f_{j-1}|^2 + \frac{2}{n} E\left(\max_{j=0,\dots,M-1} e_{i+j} - c_M\right)^2.$$

Hence,

$$(3.11) \qquad E \sum_{i=1}^{n} |\hat{f}_i - f_i|^2 \le 4M^2 \sum_{i=1}^{n-1} |f_{i+1} - f_i|^2 + 4\,\mathrm{var}\left(\max_{j=1,\dots,M} e_1\right).$$

From (3.7) and (3.8) and taking $\beta = 1$, we know that $\sum_i |f_{i+1} - f_i|^2$ is either of order $O(n^{-(2m\wedge 2)} \log n)$ or $O(n^{-(2s\wedge 2)} \log n)$, according to $p$. Thus, $\mathrm{var}(\max_{j=1,\dots,M} e_j)$ and an optimal choice of $M$ control the right-hand side of (3.11).

If the $e_i$ are i.i.d. standard normal random variables, then $\mathrm{var}(\max_{j=1,\dots,M} e_j)$ is of order $1/2 \log M$. Hence, and say, for $p \ge 2$, minimizing in $M$ ($M = n^{m\wedge 1}/(\log n)^2$) gives a rate of order $O(1/\log n)$, coming short of the thresholding rate.

Now, using arguments similar to the ones in the proof of [1], Theorem 4.1, it is easy to show that, for i.i.d. (symmetric) bounded noise, soft thresholding has the same minimax rate as it would have for Gaussian noise with the same variance. This can come short of the rate achieved by the estimator presented above. Indeed, if $e_1$ is a symmetric Bernoulli random variable with law $(\delta_{-1} + \delta_1)/2$, then $\mathrm{var}(\max_{j=1,\dots,M} e_j) = 1/2^{M-2} - 1/2^{2M-2}$. Hence, for $p \ge 2$,

$$(3.12) \qquad \sup_{\|\theta\|_{B^m_{p,q}} \le A} E \sum_{i=1}^{n} |\hat{f}_i - f_i|^2 \le C(M^2 n^{-(2m\wedge 2)} \log n + 2^{-M}),$$

where $C$ is a constant. The right-hand side in (3.12) is now minimized by choosing $M = 2\log_2 n$ and, thus,

$$(3.13) \qquad \sup_{\|\theta\|_{B^m_{p,q}} \le A} E \sum_{i=1}^{n} |\hat{f}_i - f_i|^2 = O\left(\frac{(\log n)^3}{n^{2m\wedge 2}}\right).$$

For $p \leq 2$, the right-hand side of (3.13) should be replaced by $O(\frac{(\log n)^3}{n^{2s \wedge 2}})$. In both cases the rate is better than $O(n^{-2m/(2m+1)})$, which is the minimax rate for soft thresholding in the Gaussian model.

For another example, let the $e_i$ be uniformly distributed on $[-1, 1]$. Then $\operatorname{var} \max_{j=1,\ldots,M} e_j = 4M/(M+1)^2(M+2)$ is of order $O(1/M^2)$; hence, for $p \geq 2$,

$$(3.14) \qquad \sup_{\|\theta\|_{B^m_{p,q}} \leq A} E \sum_{i=1}^n |\hat{f}_i - f_i|^2 \leq C(M^2 n^{-(2m \wedge 2)} \log n + M^{-2})$$

[resp. $\leq C(M^2 n^{-(2s \wedge 2)} \log n + M^{-2})$ for $p \leq 2$], again $C$ is a constant. For $p \geq 2$, taking $M = n^{(m \wedge 1)/2}/\sqrt[4]{\log n}$ [resp. $M = n^{(s \wedge 1)/2}/\sqrt[4]{\log n}$ for $p \leq 2$] gives

$$(3.15) \qquad \sup_{\|\theta\|_{B^m_{p,q}} \leq A} E \sum_{i=1}^n |\hat{f}_i - f_i|^2 = O\left(\frac{\sqrt{\log n}}{n^{m \wedge 1}}\right)$$

[resp. $O(\frac{\sqrt{\log n}}{n^{s \wedge 1}})$]. These rates are better [only for $p > 1/(m+1/2-2m/(2m+1))$, when $s \wedge 1 = s$] than $O(n^{-2m/(2m+1)})$. In view of [16] and of [1], Theorem 5.1, the smoothness of the density of the compactly supported noise might help thresholding reach the minimax rate among all estimators.

**4. Concluding remarks on block thresholding and kernel estimators.** Block thresholding, which applies thresholding to a whole block of wavelet coefficients, has been developed by Cai [2] as well as Hall, Kerkyacharian and Picard [11], to deal with signals exhibiting a correlation in the size of their wavelet coefficients which are above each other. More precisely, a block of noisy wavelet coefficients $\theta_1 + z_1, \ldots, \theta_k + z_k$, is kept if $\sum (\theta_i + z_i)^2$ is larger than a threshold, otherwise the whole block is set to zero (one could also keep a block if one of the coefficients in it is larger than a threshold). As defined, block thresholding shares the minimax properties of soft thresholding, in both the ideal and functional frameworks.

In block thresholding the blocks are horizontal, that is, made up of the coefficients with indices $(j, k), \ldots, (j, k + K)$. Below, we briefly present a vertical block thresholing methodology (the blocks are vertical and not disjoint) which also shares the same minimax properties as the horizontal block thresholding estimator. More importantly, we show that (for the Haar wavelet) this thresholding estimator is nothing but a kernel estimator with locally varying bandwidth. This is another instance of the well-known fact that thresholding rules represent a method of adaptive local selection of bandwidth (see [7]).

First we introduce some terminology. We say that an index $(j', k')$ (or the wavelet coefficient with this index) is above the index $(j, k)$ if $j' \leq j$ and $|[2^{j'-j}k] - k'| \leq J$, where $J \in \mathbb{N}$ is a positive constant. In vertical block thresholding, if $|\theta_{j,k} + z_{j,k}|$ is larger than a threshold, then the coefficient itself is kept and, moreover, all the coefficients above it are also kept. A variation of this method is to keep the coefficients with the indices $(j, k')$, $|k - k'| \leq J$ (for some other constant $J$), as well as all the coefficients above them.

This new method achieves (as quickly shown below) the optimal minimax rate in the ideal estimator context (a similar result holds for the function space approach too, but the proof is left out). In our usual model, let the noise be i.i.d. standard normal random variables and let $\lambda_n$ be such $E\mathbf{1}_{\{|z_{j,k}| > \lambda_n - 1\}}(1 + z_{j,k}^2) = 1/n$. With the background and methods of the present paper and its companion [1], it is easy to see that $\lambda_n \sim \sqrt{2 \log n}$. Let $J$ be as above. Note that the number of coefficients above a coefficient is less than $(2J + 1) \log_2 n$, since in each level there are only $2J + 1$ coefficients above a fixed coefficient. Next, set $Y = W_n(X)$ and define the estimator $\hat{\theta}_{j,k}$ for the coefficient $\theta_{j,k}$ by

$$\hat{\theta}_{j,k} := \begin{cases} Y_{j,k}, & \text{if } |Y_{j,k}| \geq \lambda_n \\ & \text{or } \exists (j', k'), (j,k) \text{ is above } (j', k') \text{ and } |Y_{j',k'}| \geq \lambda_n, \\ 0, & \text{elsewhere.} \end{cases}$$

We then have

$$\sum_{j,k} E(\hat{\theta}_{j,k} - \theta_{j,k})^2$$

$$\leq \sum_{j,k} E(\mathbf{1}_{\{|Y_{j,k}| \geq \lambda_n\}} z_{j,k}^2 + \mathbf{1}_{\{|Y_{j,k}| < \lambda_n\}} \theta_{j,k}^2$$

$$(4.1) \qquad\qquad + \mathbf{1}_{\{(j,k) \text{ above a } |Y_{j',k'}| \geq \lambda_n\}} z_{j,k}^2)$$

$$\leq \sum_{j,k} E\Bigg(\mathbf{1}_{\{|Y_{j,k}| \geq \lambda_n\}} z_{j,k}^2 + \mathbf{1}_{\{|Y_{j,k}| < \lambda_n\}} \theta_{j,k}^2$$

$$+ \mathbf{1}_{\{|Y_{j,k}| \geq \lambda_n\}} \sum_{(j',k') \text{ above } (j,k)} z_{j',k'}^2\Bigg).$$

If $|\theta_{j,k}| < 1$, then

$$E\mathbf{1}_{\{|\theta_{j,k} + z_{j,k}| \geq \lambda_n\}} z_{j,k}^2 \leq E\mathbf{1}_{\{|z_{j,k}| \geq \lambda_n - 1\}} z_{j,k}^2 \leq \frac{1}{n},$$

$$E\mathbf{1}_{\{|\theta_{j,k} + z_{j,k}| \leq \lambda_n\}} \theta_{j,k}^2 \leq \theta_{j,k}^2$$

and

$$E\mathbf{1}_{\{|Y_{j,k}| \geq \lambda_n\}} \sum_{(j',k') \text{ above } (j,k)} z_{j',k'}^2 \leq (2J + 1) \log_2 n E\mathbf{1}_{\{|z_{j,k}| \geq \lambda_n - 1\}} E z_{j,k}^2$$

$$\leq \frac{(2J+1)\log_2 n}{n}.$$

If $|\theta_{j,k} + z_{j,k}| < \lambda_n$, then $|\theta_{j,k}| < |\lambda_n| + |z_{j,k}|$; thus

$$E\mathbf{1}_{\{|\theta_{j,k}+z_{j,k}|<\lambda_n\}}\theta_{j,k}^2 \leq 2|\lambda_n|^2 + 2Ez_{j,k}^2.$$

Moreover,

$$E\mathbf{1}_{\{|Y_{j,k}|\geq\lambda_n\}} \sum_{(j',k') \text{ above } (j,k)} z_{j',k'}^2 \leq (2J+1)\log_2 n$$

and

$$E\mathbf{1}_{\{|\theta_{j,k}+z_{j,k}|\geq\lambda_n\}} z_{j,k}^2 \leq 1.$$

Hence

$$(4.2) \quad \frac{E(\mathbf{1}_{\{|Y_{j,k}|\geq\lambda_n\}}z_{j,k}^2 + \mathbf{1}_{\{|Y_{j,k}|<\lambda_n\}}\theta_{j,k}^2 + \mathbf{1}_{\{|Y_{j,k}|\geq\lambda_n\}} \sum_{(j',k') \text{ above } (j,k)} z_{j',k'}^2)}{1/n + \min(\theta_{j,k}^2, 1)}$$

$$\leq C\log n,$$

for some constant $C$. Combining (4.1) and (4.2), we obtain

$$\frac{\sum_{j,k} E(\hat{\theta}_{j,k} - \theta_{j,k})^2}{1 + \sum_{j,k} \min(\theta_{j,k}^2, 1)} \leq C\log n,$$

proving our claim on the minimaxity of the method.

Another interest of the vertical block thresholding method is the fact that it is close to a kernel estimate with locally varying bandwidth (this is precisely proved below in the case of the Haar wavelet). Indeed, a simple first-order approximation of the noisy wavelet coefficients is given by (since $2^j$ is small compared to $n$)

$$\tilde{\theta}_{j,k} := \sum_i \frac{\psi_{j,k}(i/n)}{\sqrt{n}} X_i,$$

where, as usual, $\psi_{j,k}$ and $\phi_{j,k}$ are, respectively, translations and dilations of the wavelet $\psi$ and of the scaling function $\phi$.

If we estimate $f_i$ by discarding the levels below the level $j_0$, then by a first-order approximation, as above,

$$\hat{f}_i := \sum_{j\geq j_0, k} \tilde{\theta}_{j,k} \frac{\psi_{j,k}(i/n)}{\sqrt{n}}$$

$$= \sum_{j\geq j_0, k} \left( \sum_l \frac{\psi_{j,k}(l/n)}{\sqrt{n}} X_l \right) \frac{\psi_{j,k}(i/n)}{\sqrt{n}}$$

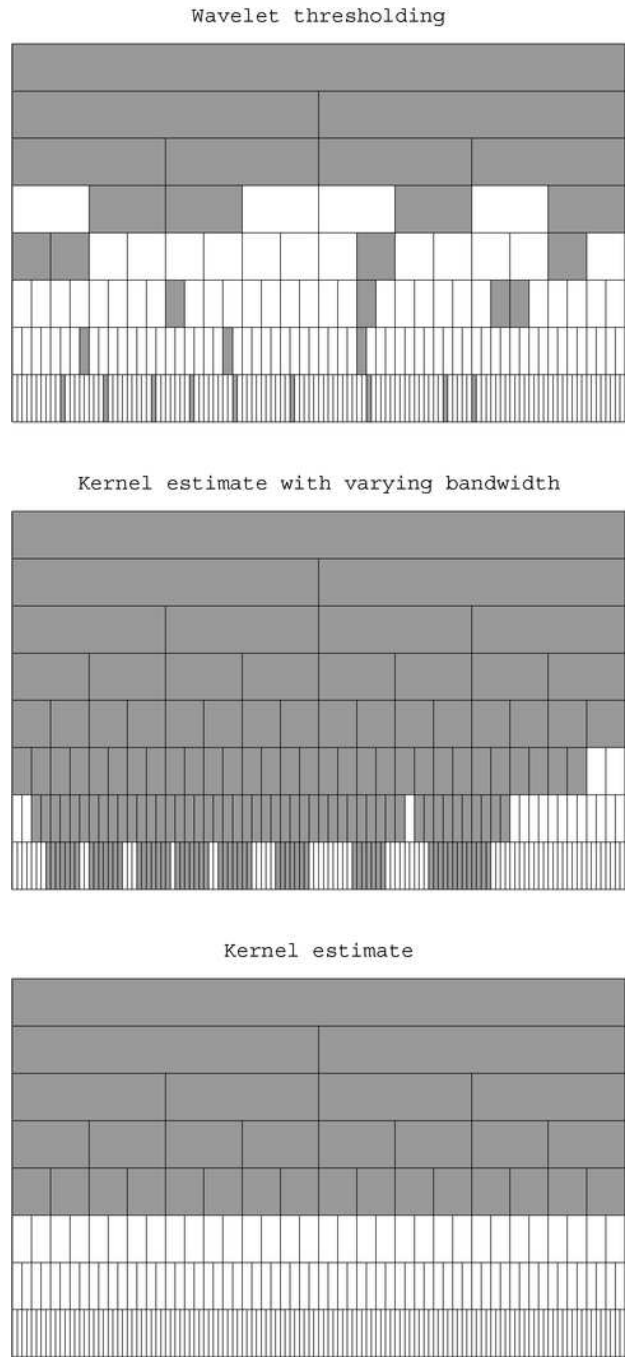

FIG. 1. *Hard thresholding, vertical block thresholding and kernel estimate.*

$$= \frac{1}{n} \sum_l X_l \sum_{j \geq j_0, k} \psi_{j,k}(l/n) \psi_{j,k}(i/n)$$

$$= \frac{1}{n} \sum_l K(l/n, i/n) X_l,$$

where $K(x,y) = \sum_{j \geq j_0, k} \psi_{j,k}(x) \psi_{j,k}(y)$. If we also keep the level $j_0 + 1$, then $K(x,y)$ has to be replaced by $K(2x, 2y)/2$. Thus, the parameter $2^{-j_0}$ corresponds to the bandwidth of a classical linear kernel estimator (see also [7]).

Figure 1 shows for an artificial signal which wavelet coefficients are kept with different methods [the artificial signal is just a random signal; the coefficient with level $j, k$ is a random variable with distribution $N(0, 2^{-\alpha j})$]. (Nothing else but these coefficients is present in the signal.) The dark rectangles correspond to coefficients which are kept. The top picture shows the coefficients kept for a hard thresholding estimator, while the bottom one shows which coefficients are kept for a kernel estimator. The middle picture illustrates why the vertical block thresholding can be viewed as a kernel estimator with locally varying bandwidth (we keep some neighboring coefficients as well).

This analogy between the vertical block thresholding estimator and kernel estimators with locally varying bandwidth becomes even more transparent by choosing the underlying wavelet basis to be the Haar basis. Then for vertical block thresholding, each estimate $\hat{f}_i$ is the mean of some neighboring $X_j$. For the Haar wavelet, the scaling identities have the forms

$$\phi_{j,k} = \frac{1}{\sqrt{2}}(\phi_{j+1,2k} + \phi_{j+1,2k+1}) \quad \text{and} \quad \psi_{j,k} = \frac{1}{\sqrt{2}}(\phi_{j+1,2k} - \phi_{j+1,2k+1}).$$

With this in mind, and for an input signal $X_0, \ldots, X_{n-1}$, $n = 2^h$, the discrete wavelet transform is given by

$$c_0 = \frac{1}{\sqrt{n}} \sum_{i=0}^{n-1} X_i \quad \text{and} \quad d_{j,k} := \frac{1}{\sqrt{2^{h-j}}} \sum_{i=0}^{2^{h-j-1}-1} f_{2^{h-j}k+i} - \sum_{i=2^{h-j-1}}^{2^{h-j}-1} f_{2^{h-j}k+i}.$$

The inverse transformation is then given by

$$f_i = \frac{1}{\sqrt{n}} c_0 + \sum_j d_{j,[i/2^{h-j}]} \begin{cases} 1, & \text{if } i/2^{h-j} - [i/2^{h-j}] < 1/2, \\ -1, & \text{if } i/2^{h-j} - [i/2^{h-j}] \geq 1/2. \end{cases}$$

To compute an estimate of $f_i$, and discarding the levels below $j_0$, we have

$$\hat{f}_i = \frac{1}{\sqrt{n}} c_0 + \sum_{j=0}^{j_0} d_{j,k} \begin{cases} 1, & \text{if } i/2^{h-j} - [i/2^{h-j}] < 1/2, \\ -1, & \text{if } i/2^{h-j} - [i/2^{h-j}] \geq 1/2, \end{cases}$$

$$= \frac{1}{2^{h-j_0-1}} \sum_{l=[i/2^{h-j_0-1}]2^{h-j_0-1}}^{[i/2^{h-j_0-1}]2^{h-j_0-1}+2^{h-j_0-1}-1} X_l,$$

where the last equality follows from a simple induction argument on $j_0$. Thus, if we discard the levels below $j_0$, then $\hat{f}_i$ is the mean of a block of $2^{h-j_0-1}$ $X_l$'s.

We claim now that, if in vertical block thresholding the coefficient with index $(j_1, [i/2^{h-j_1}])$ is kept because the coefficient with index $(j, k)$ is larger than the threshold and $|[i/2^{h-j_1}] - [k/2^{j-j_1}]| \leq J$, then for all $j_2 < j_1$, $|[i/2^{h-j_2}] - [k/2^{j-j_2}]| \leq J$, that is, the coefficients with indices $(j_2, [i/2^{h-j_2}])$, $j_2 < j_1$, are also kept.

Since for $x \in \mathbb{R}$ and $k \in \mathbb{N}$, $[x/k] = [[x]/k]$, it is clear that

$$|[i/2^{h-j_1}]/2^{j_1-j_2} - [k/2^{j-j_1}]/2^{j_1-j_2}| \leq J/2^{j_1-j_2},$$

hence

$$J \geq |[[i/2^{h-j_1}]/2^{j_1-j_2}] - [[k/2^{j-j_1}]/2^{j_1-j_2}]|$$
$$= |[i/2^{h-j_2}] - [k/2^{j-j_2}]|.$$

Thus, for vertical block thresholding, we also obtain

$$\hat{f}_i = \frac{1}{2^{h-j_0-1}} \sum_{l=[i/2^{h-j_0-1}]2^{h-j_0-1}}^{[i/2^{h-j_0-1}]2^{h-j_0-1}+2^{h-j_0-1}-1} X_l,$$

but where now $j_0$ depends on $i$ and $(X_l)$, that is, it is a kernel estimator with locally varying bandwidth.

Lepski, Mammen and Spokoiny [13] have already presented a kernel estimator with locally varying bandwidth which achieves the same minimax rate as a wavelet thresholding estimator. There the local bandwidth is chosen from a set $a^{-j}h_1$, $a, h_1 > 0$ constants and $j = 0, 1, \ldots$. For the simple kernel estimator based on wavelets, the bandwidth is $2^{-j}$, $j = 0, 1, \ldots$, and $j$ is the last level of wavelet coefficients that is kept. The results of Lepski, Mammen and Spokoiny [13] show that kernel estimates with a locally varying bandwidth selection can be as good as wavelet thresholding in a minimax sense. The performance of vertical block thresholding also makes this plausible.

## APPENDIX

Let us start with a simple lemma important in transferring part of the proof to a truncated noise setting.

LEMMA A.1.   *Let $X_i$, $i = 1, \ldots, n$, be independent random variables such that $EX_i = 0$, $EX_i^2 = 1$, and $m_4 := EX_i^4 < +\infty$. Let $Y := \sum_{i=1}^n a_i X_i$, where $\sum_{i=1}^n a_i^2 = 1$. Then $\min(3, m_4) \leq EY^4 \leq \max(3, m_4)$.*

Proof.

$$E\left(\sum_i a_i X_i\right)^4 = \sum_i a_i^4 E X_i^4 + 3 \sum_{i,j,i\neq j} a_i^2 a_j^2 E X_i^2 E X_i^2$$

$$= m_4 \sum_i a_i^4 + 3 \sum_i a_i^2 \sum_{j,j\neq i} a_j^2$$

$$= m_4 \sum_i a_i^4 + 3 \sum_i a_i^2 (1 - a_i^2)$$

$$= (m_4 - 3) \sum_i a_i^4 + 3$$

$$= m_4 + (m_4 - 3)\left(\sum_i a_i^4 - 1\right).$$

The assertion now follows from $\sum_i a_i^4 \leq 1$.  □

Proof of Theorem 2.1.  Let us describe the general strategy of proof. The risk given by the denominator of (2.3) is split into two sums going from coarser to finer noisy wavelet coefficients. First the coefficients from a certain index upward, that is, the finer wavelet coefficients, are discarded because their $\ell^2$-norm is asymptotically negligible compared to the minimax risk. Next, in view of the moment conditions imposed on $e$ and from a proper choice of the thresholds, we can reduce the proof to a truncated noise case. The rest of the proof then deals with the core of the estimation problem which corresponds to the sum containing the coarser coefficients and truncated noise. There, the right thresholds can achieve the same minimax performance as in the Gaussian case.

Choose $\alpha, \varepsilon > 0$ such that

$$\alpha > \begin{cases} 1/(2m+1), & \text{if } p \geq 2, \\ m/(2,+1)s, & \text{if } 1 \leq p < 2, \end{cases} \quad \text{and} \quad L > \frac{6}{(1-\alpha) - 2\varepsilon}.$$

This is certainly possible given the conditions of Theorem 2.1. Then let $l = l(\alpha, n)$ be such that $2^l \leq n^\alpha < 2^{l+1}$. In view of Lemma 2.2, it follows that

$$\text{(A.1)} \qquad \sup_{\|\theta\|_{B_{p,q}^m} \leq A} \sum_{j>l,k} \theta_{j,k}^2 = o(n^{-2m/(2m+1)}).$$

But [see (1.5) and the references given there] the numerator in (2.3) is $\sim C n^{-2m/(2m+1)}$. Thus, choosing $\tilde{\lambda}_{j,k} = \infty$ for $j > l$, it suffices to show that

$$\liminf_{n\to\infty} \frac{\inf_{(\lambda)\in\mathbb{R}^n} \sup_{\theta:\|\theta\|_{B_{p,q}^m}\leq A} E_\Phi \sum_{j,k}(T_{\lambda_{j,k}}^S(w_{j,k}) - \theta_{j,k})^2}{\inf_{(\tilde{\lambda})\in\mathbb{R}^{2^l}} \sup_{\theta:\|\theta\|_{B_{p,q}^m}\leq A} E \sum_{j\leq l,k}(T_{\tilde{\lambda}_{j,k}}^S(w_{j,k}) - \theta_{j,k})^2} \geq 1.$$

Let $A_n := \{\max_i |e_i| \le c_n\}$, where $c_n = n^{-\varepsilon} 2^{(h-l)/2} \le 2n^{-\varepsilon+(1-\alpha)/2}$, and let also $\tilde{e}_i := e_i \mathbf{1}_{A_n}$. Note that $\tilde{\sigma}_n^2$, the variances of the $\tilde{e}_i$, are smaller than 1, but converge to 1 (if the $e_i$ are not identically distributed, the convergence to 1 will hold uniformly). Finally, let $\tilde{z}_{j,k} := (W(\tilde{e}/\sqrt{n}))_{j,k}$.

On $A_n$, $e_i = \tilde{e}_i$, hence, denoting by $T_{j,k}$ the soft-thresholding operators with thresholds $\tilde{\lambda}_{j,k}$ smaller than $\log n/\sqrt{n}$, for $j \le l$, we have

$$E \sum_{j \le l,k} |T_{j,k}(\theta_{j,k} + z_{j,k}) - \theta_{j,k}|^2$$

$$= E\left( \sum_{j \le l,k} |T_{j,k}(\theta_{j,k} + z_{j,k}) - \theta_{j,k}|^2 \mathbf{1}_{A_n} \right.$$

$$\left. + \sum_{j \le l,k} |T_{j,k}(\theta_{j,k} + z_{j,k}) - \theta_{j,k}|^2 \mathbf{1}_{A_n^c} \right)$$

(A.2)
$$= E \sum_{j \le l,k} |T_{j,k}(\theta_{j,k} + \tilde{z}_{j,k}) - \theta_{j,k}|^2 \mathbf{1}_{A_n}$$

$$+ E \sum_{j \le l,k} |T_{j,k}(\theta_{j,k} + z_{j,k}) - \theta_{j,k}|^2 \mathbf{1}_{A_n^c}$$

$$\le E \sum_{j \le l,k} |T_{j,k}(\theta_{j,k} + \tilde{z}_{j,k}) - \theta_{j,k}|^2 \mathbf{1}_{A_n}$$

$$+ \sum_{j \le l,k} \sqrt{E|T_{j,k}(\theta_{j,k} + z_{j,k}) - \theta_{j,k}|^4 P(A_n^c)}$$

$$\le E \sum_{j \le l,k} |T_{j,k}(\theta_{j,k} + \tilde{z}_{j,k}) - \theta_{j,k}|^2$$

$$+ 2n^{\alpha}\sqrt{8(M + d_n^4)}\sqrt{P(A_n^c)},$$

where $d_n = O(\log n/\sqrt{n})$, using the elementary inequality $(a+b)^4 \le 8(a^4+b^4)$, and using Lemma A.1 [$E z_{j,k}^4 \le M := \max(3, E e_1^4)/n^2$, since $E e_i = 0$ and $E e_i^2 = 1$].

We will now show that the rightmost term in (A.2) is of order $o(1/n)$, which is again asymptotically negligible compared to the minimax risk in (2.3). Indeed, the $e_i$ have moments of order $L$; hence (using the i.i.d. assumption)

$$P(A_n^c) = P\left( \max_{1 \le i \le n} |e_i 2^{-(h-l)/2}| \ge n^{-\varepsilon} \right) \le nP(|e_1|2^{-(h-l)/2} \ge n^{-\varepsilon})$$

$$\le nE|e_1|^L n^{-L(1-\alpha)/2} n^{\varepsilon L}$$

$$\le E|e_1|^L n^{1-L((1-\alpha)/2-\varepsilon)}.$$

This implies that $P(A_n^c) = O(1/n^2)$ if $1 - L((1-\alpha)/2 - \varepsilon) \leq -2$, that is, if

$$(A.3) \qquad L \geq \frac{6}{(1-\alpha) - 2\varepsilon},$$

and this proves our claim on the size of the rightmost term in (A.2). Thus, we will be done if we prove that

$$(A.4) \quad \liminf_{n\to\infty} \frac{\inf_{(\lambda)\in\mathbb{R}^n} \sup_{\theta:\|\theta\|_{B_{p,q}^m}\leq A} E_\Phi \sum_{j,k}(T^S_{\lambda_{j,k}}(w_{j,k}) - \theta_{j,k})^2}{\inf_{\substack{(\tilde\lambda)\in\mathbb{R}^{2l} \\ \|\tilde\lambda\|_\infty \leq \log n/\sqrt{n}}} \sup_{\theta:\|\theta\|_{B_{p,q}^m}\leq A} E \sum_{j\leq l,k}(T^S_{\tilde\lambda_{j,k}}(\tilde w_{j,k}) - \theta_{j,k})^2} \geq 1.$$

Note that in (A.4) the symmetry assumption on $e$ ensures that $E\tilde e_i = 0$ for all $i$ and, thus, $E\tilde z_{j,k} = 0$ for all $j,k$.

Consider the coefficients in the levels $l$ and above with $2^l \leq n^\alpha < 2^{l+1}$. Let $\tilde z$ be the noise part in one of these coefficients. Then (see the proof of Theorem 4.1 in [1]) with $n = 2^h$,

$$(A.5) \qquad \tilde z = \sum_{i=1}^n v_i \frac{\tilde e_i}{\sqrt{n}} \quad \text{and} \quad \max_i |v_i| \leq C_1 2^{-(h-l)/2} \leq C_1\sqrt{n^{\alpha-1}},$$

where $C_1$ depends only on the type of the wavelet transform used. Note also that, by the scaling identities (e.g., see [3]),

$$(A.6) \quad \#\{v_i : v_i \neq 0\} \leq C_2 n^{1-\alpha} \quad \text{and, thus,} \quad \sum_i |v_i|^3 \leq C_1 C_2 \sqrt{n^{\alpha-1}},$$

where again $C_2$ depends only on the wavelet transform.

Since $\|\tilde e_i 2^{-(h-l)/2}\|_\infty \leq n^{-\varepsilon}$, $i = 1, \ldots, n$, the noise terms in the wavelet coefficients are sums of independent random variables which are smaller than $n^{-\varepsilon}/\sqrt{n}$, which in view of (A.5) and of (A.6) satisfy the conditions of Lemmas 2.4 and 2.5.

In the sequel, for a law $\mu$, and for $\lambda \geq 0$ and $\theta \in \mathbb{R}$, we set $p_\mu(\lambda, \theta) := \int_{-\infty}^{+\infty}(T^S_\lambda(x + \theta) - \theta)^2 \mu(dx)$. Let now $\tilde\mu_{j,k}$ denote the law of the random variable $\tilde z_{j,k}$, that is, the distribution of the noise in the coefficient of index $(j,k)$. [Recall that if the $e_i/\sqrt{n}$ are i.i.d. $N(0, 1/n)$ random variables, then the distribution of the noise in each coefficient is $\Phi_n := N(0, 1/n)$. Recall also that $Ez_{j,k}^2 = 1/n$ and, thus, $E\tilde z_{j,k}^2 \leq 1/n$, and that, finally, $l$, $\tilde\mu_{j,k}$ and $\lambda_{j,k}$ depend on $n$, but that for simplicity we choose not to indicate this in the notation.]

Let $\lambda_n$ be the threshold such that $p_{\Phi_n}(\lambda_n, 0) = 1/n^2$, and let $\theta \geq 0$. If $\lambda > \lambda_n$, then

$$(T^S_{\lambda_n}(x + \theta) - \theta)^2 < (T^S_\lambda(x + \theta) - \theta)^2 \qquad \text{for } x \in (-\lambda_n - \theta, \lambda_n).$$

Moreover, $\int_{\lambda_n}^\infty(x - \lambda_n)^2\Phi_n(dx) = p_{\Phi_n}(\lambda_n, 0)/2$ and

$$\int_{-\infty}^{-\lambda_n-\theta}((x + \lambda_n + \theta) - \theta)^2\Phi_n(dx) \leq \int_{-\infty}^{-\lambda_n}(x + \lambda_n)^2\Phi_n(dx) = p_{\Phi_n}(\lambda_n, 0)/2.$$

Hence

$$p_{\Phi_n}(\lambda_n, \theta) \leq \int_{-\lambda_n-\theta}^{\lambda_n} (T_\lambda^S(x+\theta)-\theta)^2 \Phi_n(dx) + \int_{\lambda_n}^{\infty} (T_{\lambda_n}^S(x+\theta)-\theta)^2 \Phi_n(dx)$$

$$+ \int_{-\infty}^{-\lambda_n-\theta} (T_{\lambda_n}^S(x+\theta)-\theta)^2 \Phi_n(dx)$$

$$\leq p_{\Phi_n}(\lambda, \theta) + 1/n^2.$$

From the above inequality (and a similar one for $\theta < 0$), it thus follows that

$$(A.7) \qquad \liminf_{n\to\infty} \frac{\inf_{(\lambda)\in\mathbb{R}^n} \sup_{\theta:\|\theta\|_{B_{p,q}^m}\leq A} \sum_{j,k} p_\Phi(\lambda_{j,k},\theta_{j,k})}{\inf_{\substack{(\lambda)\in\mathbb{R}^n \\ \|\lambda\|_\infty\leq\lambda_n}} \sup_{\theta:\|\theta\|_{B_{p,q}^m}\leq A} \sum_{j,k} p_\Phi(\lambda_{j,k},\theta_{j,k})} \geq 1,$$

since

$$\sum_{j\leq l,k} \frac{1}{n^2} \leq \frac{1}{n} = o(n^{-2m/(2m+1)}),$$

and again, $n^{-2m/(2m+1)}$ is the minimax rate in the Gaussian case. This shows that, without loss of generality and in the Gaussian case, we can assume that $\sup_{j\leq l,k} \lambda_{j,k} \leq \lambda_n \sim \frac{\sqrt{2\log n}}{\sqrt{n}}$. If $\tilde{z}_{j,k}$ had variance $1/n$, it would be enough in order to complete the proof of the theorem (since also $\sup_{j\leq l,k} \tilde{\lambda}_{j,k} \leq \log n/\sqrt{n}$) to show that

$$(A.8) \qquad \liminf_{n\to\infty} \inf_{j\leq l,k} \inf_{\lambda\leq\lambda_n} \sup_{\tilde{\lambda}\leq\log n/\sqrt{n}} \inf_\theta \frac{p_{\Phi_n}(\lambda,\theta)}{p_{\tilde{\mu}_{j,k}}(\tilde{\lambda},\theta)} \geq 1.$$

However, $\tilde{e}_i$ (and so $\sqrt{n}\tilde{z}_{j,k}$) has variance $\tilde{\sigma}_n^2$, which is smaller than 1 (but converges to 1) and so a further little adjustment is needed. Let $\mu_{j,k}$ be $\tilde{\mu}_{j,k}$ rescaled to have variance $1/n$. A simple differentiation under the integral shows that

$$p_{\mu_{j,k}}(\tilde{\lambda}/\tilde{\sigma}_n, \theta) \geq p_{\tilde{\mu}_{j,k}}(\tilde{\lambda}, \theta).$$

Hence, taking the sup over a larger set, in place of (A.8), it is enough to prove

$$(A.9) \qquad \liminf_{n\to\infty} \inf_{j\leq l,k} \inf_{\lambda\leq\lambda_n} \sup_{\tilde{\lambda}\leq\log n/\sqrt{n}} \inf_\theta \frac{p_{\Phi_n}(\lambda,\theta)}{p_{\mu_{j,k}}(\tilde{\lambda},\theta)} \geq 1.$$

NOTE. From now on, we set $\mu_n := \mu_{j,k}$ and also set $p_n := p_{\mu_{j,k}}$. Moreover, since performing computations with the factor $1/n$ is cumbersome, we will multiply the random variables and thresholds by $\sqrt{n}$, and the risks by $n$. The size of the fraction in (A.9) is unchanged by this transformation.

Next, we need two simple inequalities. First,

$$
\begin{aligned}
p_\Phi(\lambda, \theta) &= \theta^2 \Phi(-\lambda - \theta < x < \lambda - \theta) \\
&\quad + \int_{\lambda - \theta}^{+\infty} (x - \lambda)^2 \Phi(dx) + \int_{\lambda + \theta}^{+\infty} (x - \lambda)^2 \Phi(dx) \\
&\geq \theta^2 \Phi(-\lambda - \theta < x < \lambda - \theta) \\
&\quad + \frac{p_\Phi(\lambda, 0) + p_\Phi(\lambda + \theta \operatorname{sgn}(\theta), 0)}{2}.
\end{aligned}
\tag{A.10}
$$

For the second, let $\lambda \geq 1$. Then

$$
\begin{aligned}
p_\Phi(\lambda, 0) &\geq 2\Phi(x > \lambda + 1) \\
&\geq \frac{2}{\sqrt{2\pi}} \left( \frac{1}{\lambda + 1} - \frac{1}{(\lambda + 1)^3} \right) \exp(-(\lambda + 1)^2 / 2) \\
&\geq \frac{1}{\sqrt{2\pi}(\lambda + 1)} \exp(-(\lambda + 1)^2 / 2),
\end{aligned}
\tag{A.11}
$$

where the second inequality follows from a classic lower estimate on the standard normal distribution function (see [19], page 850).

Let us now proceed to prove (A.9). By Lemma 2.5, there exists a sequence $(\beta_n)$ converging to 1 and $\varepsilon_1 > 0$ [in view of (A.6) one can choose $\varepsilon_1 = (1 - \alpha)/2$], independent of the index of the wavelet coefficient such that, for $\alpha_n := \sqrt{\varepsilon_1 2 \log n}$, and all $c$ such that $|c| < \alpha_n$,

$$
\beta_n \leq \frac{\Phi((c, +\infty))}{\mu_n((c, +\infty))} \quad \text{and} \quad \beta_n \leq \frac{\Phi((-\infty, c))}{\mu_n((-\infty, c))}.
\tag{A.12}
$$

We distinguish two cases to prove (A.9), $\lambda < \alpha_n / 2$ and $\lambda \geq \alpha_n / 2$. Assume first that $\lambda < \alpha_n / 2$, and choose $\lambda = \bar{\lambda}$. For fixed $\lambda$, let $r_\theta(x) := (T_\lambda^S(x + \theta) - \theta)^2$. To spare us some further distinction of cases, assume that $\theta \geq 0$ (the case $\theta < 0$ leads below to similar results). Then $r_\theta$ is a function with one local minimum with value 0 at $x = \lambda$, moreover, if $\theta = 0$, then the minimum is attained at $[-\lambda, \lambda]$. Hence, $r'_\theta(x) \geq 0$ for $x \geq \lambda$ and $r'_\theta(x) \leq 0$ for $x \leq \lambda$. Thus, from

$$
\int_{-\infty}^{+\infty} r_\theta(x)\, d\Phi(x) = \int_{-\infty}^{\lambda} (-r'_\theta(x)) \Phi(x)\, dx + \int_{\lambda}^{+\infty} r'_\theta(x)(1 - \Phi(x))\, dx,
$$

and [integrating by parts with also $r_\theta(\lambda) = 0$]

$$
\begin{aligned}
\int_{-\alpha_n}^{\alpha_n} r_\theta(x)\, d\mu_n(x) &\leq \int_{-\alpha_n}^{\lambda} (-r'_\theta(x)) \mu_n((-\infty, x])\, dx \\
&\quad + \int_{\lambda}^{\alpha_n} r'_\theta(x) \mu_n([x, \infty))\, dx,
\end{aligned}
$$

and inequality (A.12), it easily follows that

$$(A.13) \qquad \frac{\int_{-\infty}^{+\infty} r_\theta(x)\, d\Phi(x)}{\int_{-\alpha_n}^{\alpha_n} r_\theta(x)\, d\mu_n(x)} \geq \beta_n.$$

Moreover, by Lemma 2.4 (with $K_n = n^{-\varepsilon}$ and $a_n = \log n$),

$$\begin{aligned}
\int_{\{|x|>\alpha_n\}} r_\theta(x)\, d\mu_n(x) &\leq \int_{\{|x|>\alpha_n\}} (\lambda + |x|)^2\, d\mu_n(x) \\
&\leq \int_{\{|x|>\alpha_n\}} 4x^2\, d\mu_n(x) \\
&\leq 4((\alpha_n^2 + 2)/c_n \exp(-\alpha_n^2 k_n/2) + o(\exp(-n^\varepsilon))) \\
&= o(p_\Phi(\lambda, 0)),
\end{aligned}$$

where $k_n = 1 - n^{-\varepsilon} \log n/2$ and where the last identity is obtained using (A.11) and $\lambda < \alpha_n/2$. Since $p_\Phi(\lambda, 0) \leq p_\Phi(\lambda, \theta)$, (A.9) holds for $\lambda < \alpha_n/2$.

Now, in the second case, $\lambda \geq \alpha_n/2$, choose the smallest $\tilde{\lambda}$ such that

$$(A.14) \qquad p_\Phi(\lambda + 1, 0) \geq p_n(\tilde{\lambda}, 0) \quad \text{and} \quad \tilde{\lambda} \geq \lambda.$$

It is a simple consequence of Lemma 2.4 and of the relation (A.11) that $\lambda/\tilde{\lambda} \to 1$ uniformly for $\lambda \geq \alpha_n/2$ (recall that we assumed $\lambda \leq \lambda_n \sim \sqrt{2\log n}$). Again we distinguish two cases.

First, let $|\theta| \leq 1$. Since $p_n(\tilde{\lambda}, \theta) \leq \theta^2 + p_n(\tilde{\lambda}, 0)$ and from (A.10), it follows that

$$\begin{aligned}
\inf_{|\theta| \leq 1} \frac{p_\Phi(\lambda, \theta)}{p_n(\tilde{\lambda}, \theta)} &\geq \inf_{|\theta| \leq 1} \frac{\theta^2 \Phi((-\lambda, \lambda)) + (p_\Phi(\lambda + 1, 0) + p_\Phi(\lambda, 0))/2}{\theta^2 + p_n(\tilde{\lambda}, 0)} \\
&\geq \Phi((-\alpha_n/2, \alpha_n/2)) \longrightarrow 1,
\end{aligned}$$

using (A.14) and $\lambda \geq \alpha_n/2$.

The case $|\theta| > 1$ is more complicated. Assume that $\theta > 1$ ($\theta < -1$ is treated in a similar fashion). Then, since $p_\Phi(\lambda, \theta) \geq p_\Phi(\lambda, 1)$, it follows that $p_\Phi(\lambda, \theta) > 1/2$ and, thus, $p_\Phi(\tilde{\lambda}, \theta) > 1/2$, if $\alpha_n$ is sufficiently large. Moreover,

$$\begin{aligned}
\int_{\{|x|>\alpha_n\}} (T_{\tilde{\lambda}}^S(x+\theta) - \theta)^2 \mu_n(dx) &\leq \int_{\{|x|>\alpha_n\}} 4(x^2 + \tilde{\lambda}^2) \mu_n(dx) \\
&= o(1),
\end{aligned}$$

since $\tilde{\lambda} \sim \lambda \leq \lambda_n \sim \sqrt{2\log n}$. As in obtaining (A.13), it is easy to see that

$$\inf_{\lambda \geq \alpha_n/2} \inf_{\theta > 1} \frac{\int_{-\infty}^{+\infty} (T_\lambda^S(x+\theta) - \theta)^2 \Phi(dx)}{\int_{-\alpha_n}^{\alpha_n} (T_\lambda^S(x+\theta) - \theta)^2 \mu_n(dx)} \longrightarrow 1,$$

thus

$$\liminf_{n\to\infty} \inf_{\lambda\geq\alpha_n/2} \inf_{\theta>1} \frac{p_\Phi(\tilde{\lambda},\theta)}{p_n(\tilde{\lambda},\theta)} \geq 1.$$

To finish the proof, and using $\tilde{\lambda}$ satisfying (A.14), we show that

$$\liminf_{n\to\infty} \inf_{\lambda\geq\alpha_n/2} \inf_{\theta>1} \frac{p_\Phi(\lambda,\theta)}{p_\Phi(\tilde{\lambda},\theta)}$$

(A.15)
$$= \liminf_{n\to\infty} \inf_{\lambda\geq\alpha_n/2} \inf_{\theta>1} \frac{\int_{-\infty}^{+\infty} (T_\lambda^S(x+\theta)-\theta)^2\Phi(dx)}{\int_{-\infty}^{+\infty} (T_{\tilde{\lambda}}^S(x+\theta)-\theta)^2\Phi(dx)}$$

$$\geq 1.$$

First, recall that

$$(T_\lambda^S(x+\theta)-\theta)^2 = \begin{cases} (x+\lambda)^2, & \text{if } x\leq -\lambda-\theta, \\ \theta^2, & \text{if } -\lambda-\theta\leq x\leq \lambda-\theta, \\ (x-\lambda)^2, & \text{if } x\geq \lambda-\theta, \end{cases}$$

and thus, if $x\leq\lambda-\theta$, then $(T_\lambda^S(x+\theta)-\theta)^2 \geq (T_{\tilde{\lambda}}^S(x+\theta)-\theta)^2$. Hence, for $\lambda>\tilde{\lambda}/2$,

$$\inf_{\theta\geq 1} \inf_{x\leq\tilde{\lambda}/2} \frac{(T_\lambda^S(x+\theta)-\theta)^2}{(T_{\tilde{\lambda}}^S(x+\theta)-\theta)^2} \geq \inf_{\theta\geq 1} \inf_{x\leq\tilde{\lambda}/2} \frac{(x-\lambda)^2}{(x-\tilde{\lambda})^2}$$

$$= \frac{(\tilde{\lambda}/2-\lambda)^2}{(\tilde{\lambda}/2-\tilde{\lambda})^2},$$

which converges to 1 since $\lambda/\tilde{\lambda}\to 1$. Finally, since

$$\int_{\tilde{\lambda}/2}^{+\infty} T_{\tilde{\lambda}}((x+\theta)-\theta)^2\Phi(dx) \leq \int_{\tilde{\lambda}/2}^{+\infty} (x-\tilde{\lambda})^2\Phi(dx)$$

$$\leq \int_{\tilde{\lambda}/2}^{+\infty} x^2\Phi(dx)$$

$$\leq \int_{\alpha_n/4}^{+\infty} x^2\Phi(dx) = o(1),$$

and since $p_\Phi(\lambda,\theta)>1/2$, the relation (A.15) holds. $\square$

**Acknowledgments.** It is a pleasure to thank a referee and an Associate Editor for valuable comments which greatly helped to improve the readability of the paper.

Institut für mathematische Stochastik
Freiburg University
Eckerstrasse 1
79104 Freiburg
Germany

Laboratoire d' Analyse et
   de Mathématiques Appliquées
CNRS UMR 8050
Université Paris XII
94010 Créteil Cedex
France
and
School of Mathematics
Georgia Institute of Technology
Atlanta, Georgia 30332
USA
e-mail: houdre@math.gatech.edu